\documentclass[]{interact}
\usepackage[plain]{algorithm}
\usepackage[noend]{algpseudocode}
\usepackage{mathtools}
\usepackage{epstopdf}
\usepackage[caption=false]{subfig}
\usepackage{amsmath,amssymb,euscript ,yfonts,psfrag,latexsym,dsfont,graphicx,bbm,color,amstext,wasysym,epsfig,subfig,parskip,textcomp}
\usepackage[numbers,sort&compress]{natbib}
\bibpunct[, ]{[}{]}{,}{n}{,}{,}% Citation support using natbib.sty
\makeatletter% @ becomes a letter
\def\NAT@def@citea{\def\@citea{\NAT@separator}}% Suppress spaces between citations using natbib.sty
\makeatother% @ becomes a symbol again

\theoremstyle{plain}% Theorem-like structures provided by amsthm.sty
\newtheorem{theorem}{Theorem}[section]
\newtheorem{lemma}[theorem]{Lemma}
\newtheorem{corollary}[theorem]{Corollary}

\theoremstyle{definition}

\theoremstyle{remark}

\graphicspath{{./},{./epsfigures/},{./figures/}}

\begin{document}

\newcommand{\tgmargin}[1]{}%{\marginpar{\color{red}\tiny\ttfamily{TG:} #1}}

\newcommand{\T}{1}
\newcommand{\R}{{\mathbb R}}
\newcommand{\mD}{{\mathbb D}}
\newcommand{\mE}{{\mathbb E}}  % this is for expectation -- we can change
\newcommand{\N}{{\mathcal N}}
\newcommand{\cR}{{\mathcal R}}
\newcommand{\cS}{{\mathcal S}}
\newcommand{\cC}{{\mathcal C}}
\newcommand{\D}{{\mathcal D}}
\newcommand{\cE}{{\mathcal E}}
\newcommand{\cF}{{\mathcal F}}
\newcommand{\cK}{{\mathcal K}}
\newcommand{\cL}{{\mathcal L}}
\newcommand{\cP}{{\mathcal P}}
\newcommand{\cQ}{{\mathcal Q}}
\newcommand{\cX}{{\mathcal X}}
\newcommand{\cH}{{\mathcal H}}
\newcommand{\cW}{{\mathcal W}}
\newcommand{\diag}{\operatorname{diag}}
\newcommand{\tr}{\operatorname{trace}}
\newcommand{\f}{{\mathfrak f}}
\newcommand{\g}{{\mathfrak g}}
\newcommand{\range}{\cR}  %{\operatorname{range}}
\newcommand{\trace}{\operatorname{trace}}
\newcommand{\argmin}{\operatorname{argmin}}

\newcommand{\ignore}[1]{}

\def\spacingset#1{\def\baselinestretch{#1}\small\normalsize}
\setlength{\parskip}{7pt}
\setlength{\parindent}{15pt}
\spacingset{1}

\newcommand{\mike}{\color{magenta}}
\newcommand{\rike}{\color{red}}
\definecolor{grey}{rgb}{0.6,0.6,0.6}
\definecolor{lightgray}{rgb}{0.97,.99,0.99}

\title
{A  variational derivation of a class of  BFGS-like methods\footnotemark[4]}

\author{
\name{Michele Pavon\textsuperscript{a}\thanks{Email: pavon@math.unipd.it}}
\affil{\textsuperscript{a}Dipartimento di Matematica ``Tullio Levi-Civita",
Universit\`a di Padova, via Trieste 63, 35121 Padova, Italy.}
}

\maketitle

\begin{abstract}
We provide a maximum entropy derivation of a new family of BFGS-like methods.  
Similar results are then derived for block BFGS methods. This also yields an independent proof of a result of Fletcher 1991 and its generalisation to the block case.
\end{abstract}

\begin{keywords}
Quasi-Newton method, BFGS method, maximum entropy problem, block BFGS.
\end{keywords}

%\begin{AMS}90C30,65K05\end{AMS}

\section{Introduction}

Suppose $f:\R^n\rightarrow\R$  is a $C^2$ function to be minimized.
Then Newton's iteration is
\begin{equation}
	x_{k+1} = x_k - [H(x_k)]^{-1} \nabla f(x_k), \quad k\in\N,
	\label{multiStep}
\end{equation}
where $H(x_k)=\nabla^2 f(x_k)$ is the Hessian of $f$ at the point $x_k$. In quasi-Newton methods, one employs instead an approximation $B_k$ of $H(x_k)$ to avoid the costly operations of computing, storing and inverting the Hessian ($B_0$ is often taken to be the identity $I_n$). These methods appear to perform well even in nonsmooth optimization, see \cite{LO}.
Instead of (\ref{multiStep}), one uses
\begin{equation}\label{quasi-step}
	x_{k+1} = x_k - \alpha_k B_k^{-1} \nabla f(x_k),\quad \alpha_k>0, \hspace{5 mm} k\in\N,
\end{equation}
 with $\alpha_k$ chosen by a line search,
 imposing the {\em secant} equation
\begin{equation}\label{secant}
y_k=B_{k+1} s_k,
\end{equation}
where 
\[y_k:=\nabla f(x_k + s_k)-\nabla f(x_k), \quad s_k:=\Delta x_k=x_{k+1}-x_k.
\]
The secant condition is motivated by the expansion
\begin{equation}
	\nabla f(x_k + s_k) \approx \nabla f(x_k) + H(x_k) s_k.
	\label{passage}
\end{equation}

For $n>1$, $B_{k+1}$ satisfying (\ref{secant}) is underdetermined. Various methods are used to find a symmetric $B_{k+1}$ that satisfies the secant equation (\ref{secant}) and is closest in some metric to the current approximation $B_{k}$. In several methods, $B_{k+1}$ or its inverse is  a rank one or two update of the previous  estimate \cite{NW}. 

Since for a strongly convex function the Hessian $H(x_k)$ is a symmetric positive definite matrix, we can think of its approximation $B_k$ as a covariance of a zero-mean, multivariate Gaussian distribution. Recall that in the case of two zero-mean multivariate normal distributions $p,q$ with nonsingular $n\times n$ covariance matrixes  $P,Q$, respectively,  the relative entropy (divergence, Kullback-Leibler index) can be derived in closed form
\[
	\mathbb{D}(p||q) = \int \log \frac{p(x)}{q(x)} p(x)dx= \frac{1}{2} \left[\log \det\left(P^{-1}Q\right) + tr{(Q^{-1}P) - n}\right].
\]
Since $P^{-1}$ and $Q^{-1}$ are the natural parameters of the Gaussian distributions, we write 
\begin{equation}
	\mathbb{D}(P^{-1}||Q^{-1}) = \frac{1}{2} \left[\log \det\left(P^{-1}Q\right) + \tr{(Q^{-1}P) - n}]\right.
\end{equation}
\section{A maximum entropy problem}\label{maxent}

Consider minimizing $\mathbb{D}(B^{-1}||B^{-1}_{k})$ over symmetric, positive definite $B$ subject to the secant equation
\begin{equation}\label{sec}
B^{-1}y_k=s_k.
\end{equation}
In \cite{F}, Fletcher indeed showed that the solution to this variational problem is provided by the BFGS iterate thereby providing a variational characterization for it alternative to Goldfarb's classical one \cite{G}, \cite[Section 6.1]{NW}. We take a different approach leading to a family of BFGS-like methods. 

First of all, observe that $B^{-1}y_k$ must be the given vector $s_k$. Thus, it seems reasonable that $B_{k+1}^{-1}$ should approximate $B_k^{-1}$ only in directions different from $y_k$. We are then led to consider the following new problem
\begin{equation}\label{newproblem}
	\min_{\{B=B^T, B>0\}} \mathbb{D}(B^{-1}||P^T_kB^{-1}_{k}P_k)
\end{equation}
subject to  (\ref{sec}), where $P_k$ is a rank $n-1$ matrix satisfying $P_ky_k=0$,
subject to the secant equation (\ref{sec}).
One possible choice for $P_k$ is the orthogonal projection 
\[P_k= I_n-\frac{y_ky_k^T}{y_k^Ty_k}=I_n-\Pi_{y_k}.
\]
Since $P_kB^{-1}_{k}P_k$ is singular, however, (\ref{newproblem}) does not make sense. Thus, to regularize the problem, we replace $P_k$ with the nonsingular, positive definite matrix $P_k^\epsilon=P_k+\epsilon I_n$.

The Lagrangian for this problem is
\begin{eqnarray}\nonumber
	\mathcal{L}(B,\lambda) = \frac{1}{2}\left[\log \det\left(B^{-1}(P_k^\epsilon)^{-1}B_{k}P_k^\epsilon\right) + tr\left(P_k^\epsilon B_{k}^{-1}P_k^\epsilon B\right) - n\right] + \lambda_k^T [B s_k -y_k]=\\\frac{1}{2}\left[\log \det\left(B^{-1}B_{k}\right)+\frac{1}{2}\log \det\left((P_k^\epsilon)^{-2}\right)+ tr\left(P_k^\epsilon B_{k}^{-1}P_k^\epsilon B\right) - n\right] + \lambda_k^T [B s_k -y_k].\nonumber
\end{eqnarray}
Observe that the term
\[\frac{1}{2}\log \det\left((P_k^\epsilon)^{-2}\right)
\]
does not depend on $B$ and  therefore plays no role in the variational analysis.
To compute the first variation of $\mathcal{L}$ in direction $\delta B$, we first recall a simple result. Consider the map $J$ defined on nonsingular, $n\times n$ matrices $M$ by $J(M)=\log|\det[M]|$ . Let  $\delta J(M;\delta M)$ denote the directional derivative of $J$ in direction $\delta M\in\R^{n\times n}$.
\
We then have the following result :
\begin{lemma}\label{chainr}{\em \cite[Lemma 2]{FP}}
If $M$ is nonsingular then, for any $\delta M\in\R^{n\times n}$,
\[
\delta J(M;\delta M)=\tr[M^{-1}\delta M].
\]
\end{lemma}
Observe also that any positive definite matrix $B$ is an interior point in the cone $\mathcal C$ of positive semidefinite matrices in any symmetric direction $\delta B\in\R^{n\times n}$. 
Imposing $\delta \mathcal{L}(B,\lambda;\delta B) = 0$ for all such $\delta B$, we get, in view of Lemma \ref{chainr},
\[\tr\left[\left(-(B^\epsilon_{k+1})^{-1}+ P_k^\epsilon B_{k}^{-1}P_k^\epsilon + 2 s_k \lambda_k^T\right)\delta B\right]=0,\quad \forall \delta B,
\]
which gives
\begin{equation}
	(B^\epsilon_{k+1})^{-1}= P_k^\epsilon B_{k}^{-1}P_k^\epsilon + 2 s_k \lambda_k^T.
\end{equation}
As $\epsilon \searrow 0$, we get the iteration
\begin{equation}\label{Biterate}
B_{k+1}^{-1}= P_kB_{k}^{-1}P_k + 2 s_k \lambda_k^T .
\end{equation}
Since $P_ky_k=0$, in order to satisfy the secant equation
\[B_{k+1}^{-1}y_k=s_k.
\] 
it suffices to choose the multiplier $\lambda_k$ so that 
\[2\lambda_k^Ty_k=1.
\]
We need, however, to also guarantee symmetry and  positive definiteness of the solution. We are then led to choose $\lambda_k$ as 
\begin{equation}\label{multiplier}
\lambda_k=\frac{s_k}{2y_k^Ts_k}.
\end{equation}
Finally, notice that, under the {\em curvature} assumption
\begin{equation}\label{curv}
y_k^Ts_k>0,
\end{equation}
if $B_k>0$,  indeed $B_{k+1}$ in (\ref{Biterate}) is symmetric, positive definite justifying the previous calculations. We have therefore established the following result.
 \begin{theorem} \label{1} Assume $B_k>0$ and $y_k^Ts_k>0$. A solution $B^*$ of 
\[\min_{\{B=B^T, B>0\}} \mathbb{D}(B^{-1}||P^T_kB^{-1}_{k}P_k), \]
 subject to constraint {\em (\ref{sec})}, in the regularized sense described above, is given by
 \begin{equation}\label{optimal}
 (B^*)^{-1}= \left(I_n-\frac{y_ky_k^T}{y_k^Ty_k}\right)B_{k}^{-1}\left(I_n-\frac{y_ky_k^T}{y_k^Ty_k}\right) +\frac{s_ks_k^T}{y_k^Ts_k}.
 \end{equation}
 \end{theorem}
 \section{BFGS-like methods}\label{BFGS-like algo}
 From Theorem \ref{1}, we get the following quasi-Newton iteration:
 \begin{eqnarray}\label{iteration1} x_{k+1} &=& x_k - \alpha_k B_k^{-1} \nabla f(x_k),\quad x_0=\bar{x},\\ B_{k+1}^{-1}&=& \left(I_n-\frac{y_ky_k^T}{y_k^Ty_k}\right)B_{k}^{-1}\left(I_n-\frac{y_ky_k^T}{y_k^Ty_k}\right) +\frac{s_ks_k^T}{y_k^Ts_k},\quad B_0=I_n.\label{iteration2}
\end{eqnarray}
Note that, for limited-memory iterations, this method has the same storage requirement as standard limited-memory BFGS, say $(s_j, y_j), j=k, k-1,\ldots, k-m+1$. Now let $v_k\in\R^n$ be any vector not orthogonal to $y_k$. Then
\begin{equation}\label{oblique}
P_k(v_k):=\frac{y_kv_k^T}{y_k^Tv_k}
\end{equation}
is an oblique projection onto $y_k$. Employing $P_k(v_k)$ and its transpose in place of $\Pi_{y_k}$  in (\ref{newproblem}) and performing the variational analysis after regularisation, we get a BFGS-like iteration 
\begin{equation}\label{OP}
B_{k+1}^{-1}= \left(I_n-P_k(v_k)\right)^TB_{k}^{-1}\left(I_n-P_k(v_k)\right) +\frac{s_ks_k^T}{y_k^Ts_k}
\end{equation}
In particular, if  $v_k=s_k$,  the corresponding oblique projection is
\[P_k(s_k)=\frac{y_ks_k^T}{y_k^Ts_k}.
\]
In such case, (\ref{OP}) is just the standard (BFGS) iteration for the inverse approximate Hessian
\begin{equation}\label{BFGS}B_{k+1}^{-1}= \left(I_n-\frac{y_ks_k^T}{y_k^Ts_k}\right)^TB_{k}^{-1}\left(I_n-\frac{y_ks_k^T}{y_k^Ts_k}\right) +\frac{s_ks_k^T}{y_k^Ts_k}.
\end{equation}
Here $T_k=I_n-P_k(s_k)$ is a rank $n-1$ matrix satisfying $T_ky_k=0$ as is $I-\Pi_{y_k}$.  We now get an alternative derivation of Fletcher's result \cite{F}.
\begin{corollary}\label{cor}Assume $B_k>0$ and $y_k^Ts_k>0$. A solution $B^*$ of 
\[\min_{\{B=B^T, B>0\}} \mathbb{D}(B^{-1}||B^{-1}_{k}), \]
 subject to constraint {\em (\ref{sec})} is given by the standard {\em (BFGS)} iteration {\em (\ref{BFGS})}.                
\end{corollary}
\proof We show that in the limit, as $\epsilon\searrow 0$,  $\mathbb{D}(B^{-1}||B^{-1}_{k})$ and  $\mathbb{D}\left(B^{-1}||\left(I_n-\frac{y_ks_k^T}{y_k^Ts_k}+ \epsilon I_n\right)^TB_{k}^{-1}\left(I_n-\frac{y_ks_k^T}{y_k^Ts_k} + \epsilon I_n\right)\right)$ only differ by terms not depending on $B$. Indeed,
\begin{eqnarray}\nonumber&&\mathbb{D}\left(B^{-1}||\left(I_n-\frac{y_ks_k^T}{y_k^Ts_k}+ \epsilon I_n\right)^TB_{k}^{-1}\left(I_n-\frac{y_ks_k^T}{y_k^Ts_k} + \epsilon I_n\right)\right)\\&&=\frac{1}{2}\left\{\log\det\left(B^{-1}B_k\right)+\log\det\left[\left(I_n-\frac{y_ks_k^T}{y_k^Ts_k} + \epsilon I_n\right)^{-1}\left(I_n-\frac{y_ks_k^T}{y_k^Ts_k} + \epsilon I_n\right)^{-T}\right]\nonumber\right.\\&&\left.+\tr\left[\left((1+\epsilon)I_n-\frac{y_ks_k^T}{y_k^Ts_k}\right)^TB_{k}^{-1}\left((1+\epsilon)I_n-\frac{y_ks_k^T}{y_k^Ts_k}\right)B\right]-n\right\}\nonumber
\end{eqnarray}
Note that, by the circulant property of the trace,
\[\tr\left[-\frac{s_ky_k^T}{y_k^Ts_k}B_k^{-1}(1+\epsilon)B\right]=\tr\left[-B\frac{s_ky_k^T}{y_k^Ts_k}B_k^{-1}(1+\epsilon)\right]
\]
It now suffices to observe that, for symmetric matrices $B$ satisfying (\ref{sec}) $Bs_k=y_k$, the products 
\[B\frac{s_ky_k^T}{y_k^Ts_k}=\frac{y_ks_k^T}{y_k^Ts_k} B=\frac{y_ky_k^T}{y_k^Ts_k}
\]
are independent of $B$.
\qed

Iterations (\ref{iteration1})-(\ref{iteration2}) and (\ref{iteration1})-(\ref{OP}) are expected to enjoy the same convergence properties as the canonical BFGS method \cite[Chapter 6]{NW}. They can, in principle, be applied also to nonsmooth cases along the lines of \cite{LO} with an exact line search to compute $\alpha_k$ at each step.
%\section{Limited memory BFGS-like methods}When $n$ is large, storing $B_k^{-1}$ is unfeasible. In limited memory BFGS, one stores instead a fixed number, say m, of values of $(\Delta x_j, y_j), j=k, k-1,\ldots, k-m+1$. Then, $\Delta x_{k} = - \alpha_k B_k^{-1} \nabla f(x_k)$ and  $B_{k+1}^{-1}$ is computed through (\ref{BFGS}) setting, for example, $B_{k-m}^{-1}=I_n$. If we do the same in iteration (\ref{optimal}), we get\begin{equation}B_{k+1}^{-1}=\end{equation}
\section{Block BFGS-like methods}
In some large dimensional problems, it is prohibitive to calculate the full gradient at each iteration. Consider for instance {\em deep neural networks}. A deep network consists of a nested composition of a linear transformation and a nonlinear one $\sigma$. In the learning phase of a deep network, one compares the predictions $y(x,\xi^i)$ for the input sample $\xi^i$ with the actual output $y^i$. This is done through a cost function $f_i(x)$, e.g.
\[f_i(x)=\|y^i-y(x;\xi^i)\|^2.
\]
The goal is to learn the {\em weights} $x$ through minimization of the empirical loss function
\[f(x)=\frac{1}{N}\sum_{i=1}^Nf_i(x).
\]
In modern datasets, $N$ can be in the millions and therefore calculation of the full gradient $\frac{1}{N}\sum_{i=1}^N\nabla f_i(x)$  at each iteration to perform gradient descent is unfeasible. One can then resort to {\em stochastic gradients} by sampling uniformly from the set $\{1,\ldots, N\}$ the index $i_k$ where to compute the gradient at iteration $k$. In alternative, one can also average the gradient over a set of randomly chosen samples called a ``mini-batch". In \cite{GG}, a so-called block BFGS was proposed.
Let $S_k$ be a {\em sketching matrix} of directions \cite{GG} and let ${\mathcal T}\subset [N]$. Rather than taking differences of random gradients, one computes the action of the sub-sampled Hessian on $S_k$ as
\[Y_k:= \frac{1}{|{\mathcal T}|}\sum_{i\in{\mathcal T}}\nabla^2f_i(x_k)S_k
\]
To update $B_k^{-1}$, we can now consider the problem
\begin{equation}\label{blockproblem}
	\min_{\{B=B^T, B>0\}} \mathbb{D}\left(B^{-1}||P_k^TB^{-1}_{k}P_k\right)
\end{equation}
where $I-P_k$ projects onto the space spanned by the columns of $Y_k$, subject to the block-secant equation
\begin{equation}\label{blocksecant}
B^{-1}Y_k=S_k.
\end{equation}
Again, one possible choice for $S_k$ is $I-\Pi_{Y_k}$ where $\Pi_{Y_k}=Y_k(Y_k^TY_k)^{-1}Y_k^T$ is the orthogonal projection. The same variational argument as in Section \ref{maxent} leads to the iteration
\begin{equation}\label{blockBFGS}
B_{k+1}^{-1}= \left(I-\Pi_{Y_k}\right)B_{k}^{-1}\left(I-\Pi_{Y_k}\right) + S_k(S_k^TY_k)^{-1}S_k^T.
\end{equation}
Another choice for $P_k$ is the oblique projection $I-Y_k(S_k^TY_k)^{-1}S_k^T$ leading to the iteration in \cite{GG} 
 \begin{equation}\label{GoldblockBFGS}
B_{k+1}^{-1}= \left(I-Y_k(S_k^TY_k)^{-1}S_k^T\right)^TB_{k}^{-1}\left(I-Y_k(S_k^TY_k)^{-1}S_k^T\right) + S_k(S_k^TY_k)^{-1}S_k^T.
\end{equation}
We then obtain a variational characterisation of the iteration (\ref{GoldblockBFGS}) alternative to the one of \cite[Appendix A]{GG} and generalizing Fletcher \cite{F}.
\begin{corollary}Assume $B_k>0$ and $S_k^TY_k>0$. A solution $B^*$ of 
\[\min_{\{B=B^T, B>0\}} \mathbb{D}(B^{-1}||B^{-1}_{k}), \]
 subject to constraint (\ref{blocksecant}) is given by $B_{k+1}$ in (\ref{GoldblockBFGS}).          
\end{corollary}
The proof is analogous to the proof of Corollary \ref{cor}.
\section{Numerical Experiments}\label{NE}
The algorithm (\ref{iteration1})-(\ref{iteration2}) has the  form:
\begin{algorithm}
    \caption{BFGS-like algorithm (\ref{iteration1})-(\ref{iteration2}) }
    \label{euclid}
    \begin{algorithmic}[1] % The number tells where the line numbering should start
        \Procedure{BFGS-like}{$f,Gf, x_0, tolerance$} 
            \State $B \gets I_d$ \Comment{$d$ is the dimension of $x_0$ and $I_d$ is the identity in $R^d$}
            \State $x \gets x_0$
            \For{$n=1,...,MaxIterations$} 
                \State $y \gets Gf(x)$
                \If{$||y||<tolerance$}
                	\State break
                \EndIf
                \State $SearchDirection \gets -By$
                \State $\alpha \gets LineSearch(f,GF,x,SearchDirection)$
                \State $\Delta x \gets \alpha \: SearchDirection$
                \State $S \gets I_d - \frac{y y^T}{y^T y}$
                \State $B \gets S^T B S + \frac{\Delta x \Delta x^T}{y^T dx}$
                \State $x \gets x + \Delta x$
            \EndFor
            \State \textbf{return} $x$
        \EndProcedure
    \end{algorithmic}
\end{algorithm}

%\paragraph{Convergence plots:\\}
\newpage
While the effectiveness of the BFGS-like algorithms introduced in Section \ref{BFGS-like algo} needs to be tested on a significant number of large scale benchmark  problems, we provide below  two examples where the BFGS-like algorithm (\ref{iteration1})-(\ref{iteration2}) appears to perform better than standard BFGS. Consider the strictly convex function $f$ on $\R^2$
\[
f(x_1,x_2) = e^{x_1-1}+e^{-x_2+1}+(x_1-x_2)^2 
\]
whose minimum point is $x^*\approx (0.8,1.2)$. Take as starting point: $(5,-7)$. Figure $1$ illustrates the decay of the error $||x^n-x^*||_2$ over $50$ iterations for the classical BFGS and for algorithm (\ref{iteration1})-(\ref{iteration2}).
\begin{figure}[H]\label{figure}
\begin{center}
\includegraphics[scale=0.5]{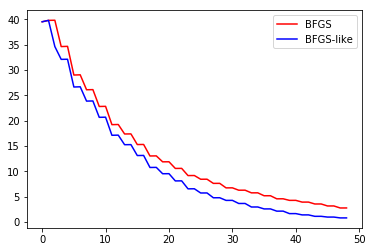}
\caption{Plot of $||x^n-x^*||_2$ for each iteration $n$}
\end{center}
\end{figure}
Consider now the (nonconvex) Generalized Rosenbrock function in $10$ dimensions:
\[f(x)=\sum_{i=1}^9\left[100\left(x_{i+1}-x_i^2\right)^2+(x_i-1)^2\right], \quad -30\le x_i\le 30, \; i=1,2,\ldots,10.
\] 
It has an absolute minumum at $x_i^*=1, i=1,\ldots,10$ and $f(x^*)=0$. Taking as initial point $x_0=(0,0,\ldots,0)$ the origin, both methods get stuck in a local minimum, see Figure $2$.
\begin{figure}[H]\label{figure2}
\begin{center}
\includegraphics[scale=0.5]{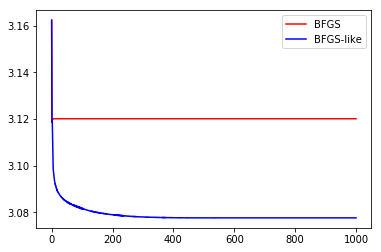}
\caption{Plot of $||x^n-x^*||_2$ for each iteration $n$}
\end{center}
\end{figure}
\newpage
Instead, initiating the recursions at $x_0=(0.9,0.9,\ldots,0.9)$, both algorithms converge to the absolute minimum (Figure $3$ depicts $100$ iterations). After a few initial steps, BFGS-like appears to perform better than BFGS.
\begin{figure}[H]\label{figure3}
\begin{center}
\includegraphics[scale=0.5]{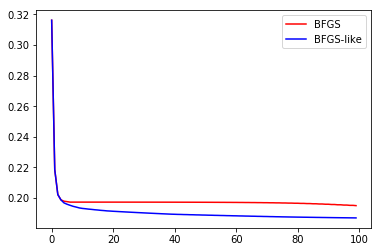}
\caption{Plot of $||x^n-x^*||_2$ for each iteration $n$}
\end{center}
\end{figure}

\section{Closing comments}
We have proposed a new family of BFGS-like iterations of which (\ref{iteration1})-(\ref{iteration2})  is a most natural one. The entropic variational derivation provides theoretical support for these methods and a new proof of Fletcher's classical derivation \cite{F}. Further study is needed to exploit the flexibility afforded by this new family (the vector $v_k$ determining the oblique projection in (\ref{oblique}) appears as a ``free parameter"). Similar results have been established for block BFGS. A few numerical experiments seem to indicate that (\ref{iteration1})-(\ref{iteration2}) may perform better in some problems than standard BFGS.

\section*{Acknowledgments}
This paper was written during a stay at the Courant Institute of Mathematical Sciences of the New York University whose hospitality is gratefully acknowledged. In particular, I would like to thank Michael Overton and Esteban Tabak for useful conversations and  for pointing out some relevant literature. I would  also like to thank Montacer Essid for kindly providing the code and the numerical examples of Section \ref{NE}.

\section*{Funding}{Supported in part by the University of Padova Research Project CPDA 140897.}


\begin{thebibliography}{99}
\bibitem{LO}A.S. Lewis and M.L. Overton,
Nonsmooth Optimization via Quasi-Newton Methods
{\em Math. Programming} {\bf 141} (2013), pp. 135-163.
\bibitem{NW} J. Nocedal and S. J. Wright, {\em Nonlinear Optimization}, 2nd edn. Springer, New York, 2006.
\bibitem{F} R. Fletcher, A New Variational Result for Quasi-Newton Formulae, {\em SIAM Journal on Optimization}, 1991, {\bf 1}, No. 1 : pp. 18-21.
\bibitem{G} D. Goldfarb, A family of variable metric methods derived by variational means, {\em Math. Comp.}, {\bf 24}, (1970), pp. 23-26.
\bibitem{FP} A. Ferrante and M. Pavon, Matrix Completion {\em \`a la}  Dempster by the Principle of Parsimony,  {\em IEEE Trans. Information Theory}, {\bf 57}, Issue 6, June 2011, 3925-3931.
\bibitem{GG}W. Gao and D. Goldfarb, Block BFGS Methods, preprint arXiv:1609.00318.

\end{thebibliography}
\end{document}